\input AHTOH-E.STY
\hfuzz 1.5pt

\UDC{
512.543.72 
+
512.542.7  
+
512.543.22 
+
512.543.27 
}
\MSC{
20F70,   
20B30,   
20E10    
}

\bigskip

\title{%
Finite symmetric groups are strongly verbally closed
}

\author{%
Olga K. Karimova$^\flat$
\qquad\qquad
Anton A. Klyachko$^{\flat\sharp}$
}
\address{
$^\flat$%
Faculty of mechanics and mathematics of Moscow State University
\\
Moscow 119991, Leninskie gory, MSU\\
$^\sharp$%
Moscow center for fundamental and applied mathematics\\
\qquad\qquad
985481@gmail.com
\qquad\qquad
klyachko@mech.math.msu.su
}

\grantsSecond{\RSF 22-11-00075}

\bigskip

\leftline{\it To the memory of Vitaly Anatol'evich Roman'kov}

\bigskip

\abstract{
\narrower
\narrower
\narrower
\narrower
Answering a question of A. V. Vasil'ev,
we show that each finite symmetric (or
alternating) group~$H$ is a retract of any
group
containing $H$ as a verbally closed subgroup.
}

\s 1.
Introduction

A subgroup $H$ of a group $G$ is called \emph{verbally closed}
[MR14]
if
any equation of the form
$
w(x,y,\dots)=h,
$
where $w$ is an element of a free group
$F(x,y,\dots)$ and $h\in H$,
having a solution in $G$ has a solution in $H$.
If each finite system of equations with coefficient from~$H$
(i.e.,
$
\{w_1(x,y,\dots)=1, \dots, w_m(x,y,\dots)=1\},
$
where $w_i\in H*F(x,y,\dots)$ and $*$ denotes the
free product)
having a solution in $G$ has a solution in $H$, then the subgroup $H$
is called \emph{algebraically closed} in $G$.

Algebraic closedness is a stronger property than verbal
closedness; but
they appear to be
equivalent surprisingly often.
The first nontrivial result on this subject
was obtained in 2014.

\proclaim Myasnikov--Roman'kov theorem {\rm[MR14]}.
Verbally closed subgroups of finite-rank free groups are
algebraically closed.

This was generalised in [KM18]:
\disp{\sl
free verbally closed subgroups of \spacing{any}
group are algebraically closed.
}%
Thus, the following notion appeared:
a group $H$ is called \emph{strongly verbally closed}
[Mazh18]
if it is algebraically closed
in any group containing $H$ as a verbally closed
subgroup.
So, verbal closedness is a property of
subgroups, while strong verbal closedness is a property of abstract
groups.

The class of strongly verbally closed groups turns out to be fairly wide.
For example, the following groups are strongly verbally closed:
\-
all abelian groups [Mazh18],
\-
all free groups [KM18],
\itemitem{--}
and even all virtually free groups
containing no
non-identity finite normal subgroups
[KM18], [KMM18]
(recall that a \emph{virtually free} group
is a group containing a free
finite-index subgroup),
\-
all free products of nontrivial groups [Mazh19],
\-
the fundamental groups of all connected surfaces,
except for the Klein bottle [Mazh18], [K21],
\itemitem{--}
and even all acylindrically hyperbolic groups
without non-identity finite normal subgroups [Bog22]
(this generalises
several results above),
\-
all finite groups with nonabelian monoliths (in particular, all
finite simple groups) [KMO23],
\-
all dihedral groups whose orders are not divisible by eight
[KMO23].

\enditem
There is also a general embedding theorem [KMO23]:
\disp{\sl
any group $H$ embeds into a strongly verbally closed group
of
cardinality~$|H|+\aleph_0$ that
satisfies all identities of $H$.
}%
The class of all groups satisfying all identities (laws) of a group $H$
is called the \emph{variety generated by $H$} and denoted
$\var H$ (see [Neu69]).
The role of identities is explained
by the following (easy-to-prove)
remark [MR14]:
\disp{\sl
if a group $H$ is a retract of any
finitely generated
group $G\in\var H$
containing $H$ as a verbally closed subgroup,
then $H$ is strongly verbally closed.
}%
(Recall that a \emph{retract}
of a group is the image of an endomorphism
$\rho$ such that $\rho\o\rho=\rho$.)
This sufficient condition of strong verbal closedness seems too
restrictive, but actually, many groups have even
stronger property: a group $H$ is called a \emph{strong retract} [KMO23]
if it is a retract of \spacing{any} group $G\supseteq H$
from~$\var H$.

\Th KMO {\rm [KMO23]}.
The following groups are strong retracts:
\-
all finite groups with nonabelian monoliths
\(in particular, all finite simple groups\)
\-
and each finite group $H$
containing a normal subgroup $C$ such that
$C$
coincides with its centraliser, does not decompose
into a direct
product of
nontrivial
normal in $H$ subgroups, and $\GCD(|C|,|H/C|)=1$.

Recall that the \emph{monolith} of a group
is the intersection of all its non-identity normal subgroups
(and a group with a non-identity monolith is called \emph{monolithic}\/).

Surely, strong retractness
is much stronger than strong verbal closedness.
For instance, all abelian groups are strongly verbally closed
(as was mentioned above),
but only few of them are strong retracts, as
the following theorem shows.

\proclaim Nilpotent-strong-retract theorem {\rm [D24]}.
A nilpotent group is a strong retract if and only if
it is abelian and
\-
either divisible,
\-
or decomposes into a direct sum of cyclic groups
of uniformly bounded
orders, which are pairwise
equal-or-coprime.

This paper answers a question of
Andrey Viktorovich Vasil'ev:
are all
finite symmetric groups strongly verbally closed?
For a majority of such groups, the answer follows from Theorem KMO.
The case of the symmetric group of degree four seems
most difficult.
However, it turns out that
some well-known (but forgotten) facts imply easily
that this group is not an exception (see the following section).

\proclaim Main theorem.
All finite symmetric and alternating groups
are strong retracts
\(in particular, they are strongly verbally closed\).

\Question.
Are infinite symmetric \(finitary and full\)
and alternating groups strongly verbally closed?

An inexperienced reader could get the impression that
(almost)
all groups are strongly verbally closed. Indeed, it is
not easy to prove the absence of this
property for any particular group.
Until quite lately, only few examples of
non-strongly-verbally-closed groups
were known:
\-
two non-abelian groups of order eight
[KM18], [RKhK17]
\-
and the fundamental group of the Klein bottle
$\pres<a,b|a^2=b^2>$ [K21].

\enditem
Strict necessary conditions for strong verbal closedness were obtained
only recently:
\disp{\sl\hfuzz11pt
the centre of a finitely generated strongly verbally closed group is
\-
a direct factor if the group is finite {\rm[KMO23]}\;
\-
pure \(= servant\) if the group is, e.g., linear
\newline
\rm (and generally ``practically always", i.e., under
some additional assumptions which are satisfied almost always)
[DenK24].
}%
For example, the following groups are not strongly verbally closed:
\-
all non-abelian nilpotent finite
groups
(and even all finitely generated
nilpotent groups with non-abelian periodic parts [D24]);
\-
all non-abelian braid groups;
\-
the group $\SL_{2024}(\Z)$.

\enditem
Other results on the verbal closedness can be found,
e.g., in
[Rom12],
[RKh13],
[Mazh17],
[RT20],
and
[Tim21].

The authors thank A. Yu. Olshanskii for a valuable remark.
The second author thanks
the Theoretical Physics and Mathematics Advancement Foundation ``BASIS".

\s 2.
Proof of the main theorem

We call a group $H$ \emph{maximal monolithic} in
a class of groups ${\cal K}\ni H$ if does not embed
into a larger group $G\in\cal K$ in such a way that
the monolith of $H$ embeds into the monolith of $G$
(i.e., the inclusions $H\subseteq G\in\cal K$
and $(\hbox{monolith of $H$)}\subseteq(\hbox{monolith of $G$})$
imply that $G=H$).

\proclaim Monolithic-strong-retract lemma.
A finite monolithic group $H$ is a strong retract
if and only if
it is maximal monolithic in the class
of finite groups from $\var H$.

\Proof
The ``only if" assertion is obvious:
the kernel of a retraction $G\to H$ onto a proper subgroup must
contain the monolith of $G$, therefore,
the monolith of $G$ must trivially intersect $H$.

Let us prove the ``if" part. Suppose that a group $H$ is maximal
monolithic in the class of finite groups from
$\var H\ni G\supseteq H$.
We have to construct a retraction~$G\to H$.
Let us apply the following
Proposition from Section 2 of~[KMO23]:
\disp{\narrower\sl
a finite subgroup $A$ of a group $B$
is a retract if and only if it is a retract
of each finitely generated subgroup of $B$ containing $A$.
}%
Thus, the group $G$ can be assumed finitely generated,
and even finite, because the variety generated by a finite group
consists of locally finite groups [Neu69].
Let us choose a maximal normal subgroup $N\nin G$ trivially
intersecting $H$ (or, equivalently, not containing the monolith $M$ of
$H$).
Then
\-
the natural homomorphism $\pi\:G\to G/N$ is injective on $H$;
\-
the group $G/N$ is monolithic with monolith containing $\pi(M)$
(by virtue of the maximality of $N$);
\-
hence, $G/N=\pi(H)$ (because $\pi(H)\iso H$ is maximal
monolithic in the class of finite groups from~$\var H$);
\-
therefore, the composition
$G\too^\pi G/N=\pi(H)\iso H$
(with a suitable isomorphism) is a retraction~$G\to H$.

\enditem
This completes the proof of the lemma.

Let us prove the main theorem.
Symmetric groups of degree at most two
and alternating groups of degree at most three
are cyclic and, therefore, strong retracts by
Denissov's nilpotent-strong-retract theorem (see Introduction).

Symmetric and alternating groups of
finite degree $n\ge5$ are monolithic with nonabelian monolith,
therefore they are strong retracts by Theorem KMO.

The monoliths of the symmetric group $S_3$ of degree three and alternating
group $A_4$ of
degree four are abelian (they are the alternating group $A_3$ of degree
three and the Klein four-group $V_4$). So, the first assertion of
Theorem~KMO does not apply, but the second assertion do apply
(take $C$ to be
the monolith).

The only remaining case is the symmetric group
$H=S_4$ of degree four.
This group is also monolithic (the monolith is the Klein four-group
$V_4$). The classification of monolithic groups from the variety
generated by the symmetric group of degree four is known [COP70]:

\dispno{\sl\hfuzz 5.4pt
finite non-nilpotent monolithic groups from $\var S_4$ are precisely
the symmetric groups
of degrees three and four and the alternating group of degree four.
}(*)%
Thus,
$S_4$ is a maximal monolithic group in the class of finite groups
from $\var S_4$. Applying the monolithic-strong-retract lemma, we
complete the proof.

\medskip

\noindent{\bf Some explanations on
fact $(*)$.}
In [COP70], a more general result was proven.

\Lemma {\rm([COP70], Lemma 4.3.1 = \emph{``an important lemma"})}.
If a non-nilpotent finite monolithic group
$G\in \eufm A_2\eufm A_3\eufm A_2$, whose
Sylow 2-subgroup lies in $\eufm D$, satisfies the
identity $\[[x,y]^3,y^3,y^2\]=1$, then $G$ is isomorphic to $S_4$, $A_4$,
or $S_3$.
\newline
{\rm Here $\eufm A_n$ is the variety of abelian groups
of exponent $n$, and $\eufm D$ is the variety generated by the dihedral
group of order eight
(this is the same as the variety generated by the quaternion group,
or the variety defined by laws
$x^4=1$ and $[x^2,y]=1$).}

\-
Surely, the variety generated by $S_4$ is contained in
the product of varieties $\eufm A_2$, $\eufm A_3$, and~$\eufm A_2$,
because the group $S_4$ itself lies in this product: the factors
of the subnormal series $\1\nin V_4\nin A_4\nin S_4$
are abelian groups of exponents two, three, and two.

\-
The symmetric group $S_4$
(and, therefore, any group from $\var S_4$)
certainly satisfies the law
$\[[x,y]^3,y^3,y^2\]=1$ (this can be easily
verified; actually, [COP70] contains an explicit basis of laws of
$S_4$).

\-
The Sylow 2-subgroup of any
finite group from $\var S_4$ lies in $\eufm D$. 
Indeed,
the Sylow 2-subgroup of $S_4$ is the dihedral group of order eight;
and there is a general fact (an immediate corollary of 
the Birkhoff variety theorem): 
\disp{\sl
the Sylow $p$-subgroup of any finite group from the variety generated by a 
finite group $G$ belongs to the variety generated by the Sylow 
$p$-subgroup of $G$.  
}

\enditem
Thus, we see that $(*)$ is indeed a
special case of the ``important lemma" from~[COP70].


\References

[Bog22]
O. Bogopolski,
Equations in acylindrically hyperbolic groups and verbal closedness,
Groups, Geometry, and Dynamics, 16:2 (2022), 613-682.
\arXiv 1805.08071

[COP70]
J. Cossey, S. Oates MacDonald, and A. Penfold Street,
On the laws of certain finite groups,
Journal of the Australian Mathematical Society, 11:4 (1970), 441-489.

[D24]
F. D. Denissov,
Finite normal subgroups of strongly verbally closed groups,
Journal of Group Theory, 27:5 (2024), 1039-1057.
\arXiv:2301.02752

[DenK24]
F. D. Denissov, A. A. Klyachko
The centre of a finitely generated strongly verbally closed group
is almost always pure,
Quarterly Journal of Mathematics, 75:3 (2024), 1149-1156.
\arXiv:2308.06837

[K21]
A. A. Klyachko,
The Klein bottle group is not strongly verbally closed,
though awfully close to being so,
Canadian Mathematical Bulletin, 64:2 (2021), 491-497.
\arXiv 2006.15523

[KM18]
A. A. Klyachko, A. M. Mazhuga,
Verbally closed virtually free subgroups,
Sb. Math., 209:6 (2018), 850-856.
\arXiv 1702.07761

[KMM18]
A. A. Klyachko, A. M. Mazhuga, V. Yu. Miroshnichenko,
Virtually free finite-normal-subgroup-free groups
are strongly verbally closed,
J. Algebra, 510 (2018), 319-330.
\arXiv 1712.03406

[KMO23]
A. A. Klyachko, V. Yu. Miroshnichenko, A. Yu. Olshanskii,
Finite and nilpotent strongly verbally closed groups,
Journal of Algebra and Its Applications 22:09 (2023), 2350188.
\arXiv:2109.12397

[Mazh17]
A. M. Mazhuga,
On free decompositions of verbally closed subgroups
of free products of finite groups,
J.~Group Theory, 20:5 (2017), 971-986.
\arXiv 1605.01766

[Mazh18]
A. M. Mazhuga,
Strongly verbally closed groups,
J. Algebra, 493 (2018), 171-184.
\arXiv 1707.02464

[Mazh19]
A. M. Mazhuga,
Free products of groups are strongly verbally closed,
Sb. Math., 210:10 (2019), 1456-1492.
\arXiv 1803.10634

[MR14]
A. Myasnikov, V. Roman'kov,
Verbally closed subgroups of free groups,
J. Group Theory, 17:1 (2014), 29-40.
\arXiv 1201.0497

[Neu69]
H. Neumann,
Varieties of groups,
Springer-Verlag, Berlin-Heidelberg-New York, 1967.

[Rom12]
V. A. Roman'kov,
Equations over groups,
Groups - Complexity - Cryptology,
4:2 (2012), 191-239.

[RKh13]
V. A. Roman'kov, N. G. Khisamiev.
Verbally and existentially closed subgroups of free nilpotent groups.
Algebra and Logic, 52:4 (2013), 336-351.

[RKhK17]
V. A. Roman'kov, N. G. Khisamiev, A. A. Konyrkhanova,
Algebraically and verbally closed subgroups and retracts
of finitely generated nilpotent groups,
Siberian Math. J., 58:3 (2017), 536-545.

[RT20]
V. A. Roman'kov, E. I. Timoshenko,
Verbally Closed Subgroups of Free Solvable Groups,
Algebra and Logic, 59:3 (2020), 253-265.
\arXiv 1906.11689

[Tim21]
E. I. Timoshenko,
Retracts and verbally closed subgroups
with respect to relatively free soluble groups,
Sib. Math. J., 62:3, 537-544 (2021).

\end